\def\Z{\mathbb Z}
\def\R{\mathbb R}
\def\Q{\mathbb Q}
\def\N{\mathbb N}
\def\A{\mathcal A}

\def\C{\mathcal C}

\def\pf{\begin{proof}}
\def\pfk{\end{proof}}

\documentclass[a4paper,11pt,reqno,oneside]{article}

\usepackage{amsmath}
\usepackage{amsthm}
\usepackage{amssymb}
\usepackage{enumerate}
\usepackage{epic,eepic}

\newtheorem{lem}{Lemma}[section]
\newtheorem{prop}[lem]{Proposition}
\newtheorem{coro}[lem]{Corollary}
\newtheorem{thm}[lem]{Theorem}
\newtheorem{de}[lem]{Definition}
\newtheorem{pozn}[lem]{Remark}

\begin{document}

\begin{center}

{\LARGE\textbf{Characterization of substitution invariant 3iet
words}}

\vskip0.5cm

{\large Peter Bal\'a\v zi, \ Zuzana~Mas\'akov\'a, \
Edita~Pelantov\'a}

\vskip0.2cm
{\small Doppler Institute 
\& Department of Mathematics, FNSPE}\\
{\small Czech Technical University, Trojanova 13, 120 00 Praha 2, Czech Republic}\\
{\small E-mails: \texttt{peter\underline{ }balazi@centrum.cz}, \
\texttt{masakova@fjfi.cvut.cz}}, \ {\small
\texttt{pelantova@fjfi.cvut.cz}}
\date{\today}

\vskip0.5cm

\begin{abstract}
We study infinite words coding an orbit under an exchange of three
intervals which have full complexity $\C(n)=2n+1$ for all $n\in\N$
(non-degenerate 3iet words). In terms of parameters of the
interval exchange and the starting point of the orbit we
characterize those 3iet words which are invariant under a
primitive substitution. Thus, we generalize the result recently
obtained for sturmian words.
\end{abstract}

\end{center}

\section{Introduction}

 We study invariance under substitution of infinite
words coding exchange of three intervals with permutation (321).
These words, which are here called 3iet words, are one of the
possible generalizations of sturmian words to a three-letter
alphabet. Our main result provides necessary and sufficient
conditions on the parameters of a 3iet word to be invariant under
substitution.

A sturmian word $(u_n)_{n\in\N}$ over the alphabet $\{0,1\}$ is
defined as
$$
u_n = \lfloor (n+1)\alpha + x_0\rfloor - \lfloor
n\alpha+x_0\rfloor \qquad \text{for all $n\in\N$,}
$$
 or
$$
u_n = \lceil (n+1)\alpha + x_0\rceil - \lceil n\alpha+x_0\rceil
\qquad \text{for all $n\in\N$,}
$$
where $\alpha\in(0,1)$ is an irrational number called the slope,
and $x_0\in[0,1)$ is called the intercept.

There are many various equivalent definitions of sturmian words,
among others also as an infinite word coding an exchange of 2
intervals of length $\alpha$ and $1-\alpha$. A direct
generalization of this definition are infinite words coding
exchange of $k$ intervals, as introduced by Stepin and
Katok~\cite{SteKat}.

\begin{de}
Let $\alpha_1,\dots,\alpha_k$ be positive real numbers and let
$\pi$ be a permutation over the set $\{1,2,\dots,k\}$. Denote
$I=I_1\cup I_2 \cup \cdots\cup I_k $, where
$I_j:=\bigl[\sum_{i<j}\alpha_i,\sum_{i\leq j}\alpha_i\bigr)$. Put
$t_j:= \sum_{\pi(i)<\pi(j)}\alpha_i - \sum_{i<j}\alpha_i$. The
mapping $T:I\mapsto I$ given by the prescription
$$
T(x)=x+t_j \quad \hbox{for } x\in I_j
$$
will be called $k$-interval exchange transformation ($k$-iet) with
permutation $\pi$ and parameters $\alpha_1,\dots,\alpha_k$.
\end{de}


Keane~\cite{kean} has studied in which case a $k$-iet satisfy the
so-called minimality condition, i.e., when the orbit
$\{T^n(x_0)\mid n\in\Z\}$ of every point $x_0\in I$ is dense in
$I$. It is easy to see that minimality condition can be satisfied
only if the permutation $\pi$ is irreducible, i.e.,
$$
\pi\{1,2,\dots,j\}\neq \{1,2,\dots,j\} \quad \hbox{for all }
j<k\,.
$$
Keane has also derived a sufficient condition for minimality:
Denote $\beta_j$ the left end-point of the interval $I_j$, i.e.,
$\beta_j=\sum_{i<j}\alpha_i$. If the orbits of points
$\beta_1,\dots,\beta_k$ under the transformation $T$ are infinite
and disjoint, then $T$ satisfies the minimality property. In the
literature, this sufficient condition is known under the notation
i.d.o.c. However, in general, i.d.o.c.\ is not a necessary
condition for minimality property.

To the orbit of every point $x_0\in I$, one can naturally
associate an infinite word $u=(u_n)_{n\in\Z}$ in a $k$-letter
alphabet $\A=\{1,2,\dots,k\}$. For $n\in\Z$ put
$$
u_n=i \quad \hbox{if } \ T^n(x_0)\in I_i\,.
$$
Infinite words coding $k$-iet with i.d.o.c.\ are called here {\em
non-degenerate $k$-iet words}. Non-degene\-rate $k$-iet words are
studied in~\cite{Ferenczi}. The authors give a combinatorial
characterization of the language of infinite words which
correspond to $k$-iet with permutation
\begin{equation}\label{eq:cislo1}
\pi(1)=k,\ \pi(2)=k-1, \ \dots\ ,\ \pi(k)=1
\end{equation}
or to permutations in some sense equivalent with it.

For $k=2$, the only irreducible permutation is of the
form~\eqref{eq:cislo1}. The minimality property for parameters
$\alpha_1,\alpha_2$ means that they are linearly independent over
$\Q$. Infinite words coding 2iet with the minimality property are
precisely the sturmian words.

In this paper we concentrate on infinite words coding exchange of
3 intervals under the permutation~\eqref{eq:cislo1}. The
transformation which we study is thus given by a triple of
positive parameters $\alpha_1,\alpha_2,\alpha_3$ and the
prescription
\begin{equation}\label{eq:cislo4}
\widetilde{T}(x):=\left\{
\begin{array}{ll}
x+\alpha_2+\alpha_3 &\hbox{for }\ x\in[0,\alpha_1)\,,\\
x-\alpha_1+\alpha_3 &\hbox{for }\ x\in[\alpha_1,\alpha_1+\alpha_2)\,,\\
x-\alpha_1-\alpha_2 &\hbox{for }\
x\in[\alpha_1+\alpha_2,\alpha_1+\alpha_2+\alpha_3)\,.
\end{array}
\right.
\end{equation}
For such a transformation, the minimality property is equivalent
to the following condition (as proved in~\cite{3iet}): numbers
$\alpha_1+\alpha_2$ and $\alpha_2+\alpha_3$ are linearly
independent over $\Q$.\footnote{Let us mention that the question
of expressing the minimality property in terms of parameters
$\alpha_1,\dots,\alpha_k$ has not been solved for general $k$.} It
is known~\cite{adamczewski,GuMaPe} that infinite words
coding~\eqref{eq:cislo4} are non-degenerate if and only
if~\eqref{eq:cislo4} satisfies the minimality property and
\begin{equation}\label{eq:cisloCC}
\alpha_1+\alpha_2+\alpha_3 \notin (\alpha_1+\alpha_2)\Z +
(\alpha_2+\alpha_3)\Z\,.
\end{equation}

The central problem of this paper is the substitution invariance
of given infinite words. For sturmian words this question was
extensively studied; Chapter X.\ of \cite{lothaire2} gives
references to authors who gave some contributions to its solution.
Complete answer to this question was first provided by
Yasutomi~\cite{Yasutomi}, other proofs of the same result are
given in~\cite{BaMaPe,berthe}. Crucial for stating this result is
the notion of a Sturm number. The original definition of a Sturm
number used continued fractions. In 1998, Allauzen~\cite{allauzen}
has provided a simple characterization of Sturm numbers:\\

{\em A quadratic irrational number $\alpha$ with conjugate
$\alpha'$ is called a Sturm number if
$$
\alpha\in(0,1) \quad\hbox{ and }\quad\alpha'\notin(0,1)\,.
$$
}

\begin{thm}[\cite{BaMaPe}]\label{t:sturmian}
Let $\alpha$ be an irrational number, $\alpha\in(0,1)$,
$x_0\in[0,1)$. The Sturmian word with slope $\alpha$ and intercept
$x_0$ is invariant under a substitution if and only if the
following three conditions are satisfied:
\begin{itemize}

\item[{\rm (i)}]
$\alpha$ is a Sturm number,

\item[{\rm (ii)}]
$x_0\in\Q(\alpha)$,

\item[{\rm (iii)}]
$\min(\alpha',1-\alpha')\leq x'_0 \leq \max(\alpha',1-\alpha')$,
where $x'_0$ denotes the image of $x_0$ under the Galois
automorphism of the quadratic field $\Q(\alpha)$.
\end{itemize}
\end{thm}

Let us mention that one can study also weaker property than
substitution invariance; namely the substitutivity. For an
infinite word $u$ coding an exchange of $k$ intervals,
Boshernitzan and Carol~\cite{Boshernitzan} have shown that
belonging of lengths of all intervals $I_1,\dots,I_k$ to the same
quadratic field is a sufficient condition for substitutivity of
$u$. For $k=2$ in~\cite{berthe2}, and for $k=3$
in~\cite{adamczewski} the respective authors show that such
condition is also necessary.

However, quadraticity of parameters is not sufficient for the
property of substitution invariance. Already in~\cite{bluetheorem}
it is shown that substitution invariance of 3iet words implies
that a certain parameter of the 3iet is a Sturm number, namely
$$
\varepsilon =
\frac{\alpha_1+\alpha_2}{\alpha_1+2\alpha_2+\alpha_3}\,.
$$
The main result of this paper is given as Theorem~\ref{t:hlavni},
where a necessary and sufficient condition for substitution
invariance is expressed using simple inequalities for other
parameters of the 3iet word.

\section{Basic notions of combinatorics on words}

We will deal with infinite words over a finite alphabet, say
$\A=\{1,2,\dots,k\}$. We consider either right sided infinite
words
$$
u=(u_n)_{n\in\N}=u_0u_1u_2u_3\cdots\,,\quad u_i\in\A\,,
$$
or {\em pointed bidirectional infinite words},
$$
u=(u_n)_{n\in\Z}=\cdots u_{-2}u_{-1}|u_0u_1u_2u_3\cdots\,,\quad
u_i\in\A\,,
$$
A finite word $w=w_0w_1\cdots w_{n-1}$ of length $|w|=n$ is a {\em
factor} of an infinite word $u=(u_n)$ if $w=u_iu_{i+1}\cdots
u_{i+n-1}$ for some $i$.

The (factor) complexity of $u=(u_n)_{n\in\N}$ is the function
$\C:\N\mapsto\N$,
$$
\C(n):=\#\{u_i\cdots u_{i+n-1}\mid i\in\N\}\,,
$$
analogously for $u=(u_n)_{n\in\Z}$. Obviously, every infinite word
satisfies $1\leq \C(n) \leq k^n$ for all $n\in\N$. It is not
difficult to show~\cite{Hedlund} that an infinite word
$u=(u_n)_{n\in\N}$ is eventually periodic if and only if there
exists $n_0$ such that $\C(n_0)\leq n_0$. Obviously, the aperiodic
words of minimal complexity satisfy $\C(n)=n+1$ for all $n\in\N$.
Such infinite words are called Sturmian words. The definition of
Sturmian words can be extended also to bidirectional infinite
words $(u_n)_{n\in\Z}$, requiring except of $\C(n)=n+1$ for all
$n\in\N$ also the irrationality of the densities of letters.

In our paper we study invariance of infinite words under
substitution. A substitution is a mapping
$\varphi:\A^*\mapsto\A^*$, where $\A^*$ is the monoid of all
finite words including the empty word, satisfying
$\varphi(vw)=\varphi(v)\varphi(w)$ for all $v,w\in\A^*$. In fact,
a substitution is a special case of a morphism
$\A^*\mapsto{\mathcal B}^*$, where $\A={\mathcal B}$. Obviously,
$\varphi$ is uniquely determined, if defined on all the letters of
the alphabet. A substitution $\varphi$ is called primitive, if
there exists $n\in\N$ such that $\varphi^n(a)$ contains $b$ for
all letters $a,b\in\A$.

The action of $\varphi$ can be naturally extended to infinite
words. For a pointed bidirectional infinite word
$u=(u_n)_{n\in\Z}$ we in particular have
$$
\varphi(\cdots u_{-2}u_{-1}|u_0u_1u_2\cdots) \ = \
\cdots\varphi(u_{-2})\varphi(u_{-1})|\varphi(u_{0})\varphi(u_{1})\varphi(u_{2})\cdots
$$
An infinite word $u$ is said to be a fixed point of $\varphi$ (or
invariant under $\varphi$), if $\varphi(u)=u$.

\section{Exchange of three intervals and cut-and-project sets}

Our aim is to study substitution invariance of words coding an
exchange of three intervals~\eqref{eq:cislo4}. The main tool is
the fact that the orbit of an arbitrary point under this
transformation can be geometrically represented by a so-called
cut-and-project sequence.

\begin{de}
Let $\varepsilon,\eta\in\R$, $\varepsilon\neq -\eta$,
$\varepsilon,\eta$ irrational, and let $\Omega=[c,c+l)$, $c\in\R$,
$l>0$. The set
\begin{equation}\label{eq:cislo2}
\Sigma_{\varepsilon,\eta}(\Omega):=\{a+b\eta\mid a,b\in\Z,\
a-b\varepsilon\in\Omega\}
\end{equation}
is called a cut-and-project set with parameters $\varepsilon,\eta$
and acceptance window $\Omega$.
\end{de}

The above definition is a very special case of a general
cut-and-project set, introduced in~\cite{moody}. In the definition
we have used an interval $\Omega$, closed from the left and open
from the right. One can also consider an interval
$\hat{\Omega}=(\hat{c},\hat{c}+\hat{l}]$. However, by doing this,
we do not obtain anything new, since
$\Sigma_{\varepsilon,\eta}(\Omega)=-\Sigma_{\varepsilon,\eta}(-\hat{\Omega})$.

For simplicity of notation, we denote the additive group
$\{a+b\varepsilon\mid
a,b\in\Z\}=\Z+\varepsilon\Z=:\Z[\varepsilon]$ and analogously for
$\Z[\eta]$. The morphism of these groups $\ x=a+b\eta \ \mapsto\
x^*=a-b\varepsilon\ $ will be called the star map. In this
formalism, the cut-and-project set
$\Sigma_{\varepsilon,\eta}(\Omega)$ can be rewritten as
$$
\Sigma_{\varepsilon,\eta}(\Omega)=\{x\in\Z[\eta]\mid
x^*\in\Omega\}\,.
$$

The relation between the set $\Sigma_{\varepsilon,\eta}(\Omega)$
and the exchange of 3 intervals is explained by the following
theorem proved in~\cite{GuMaPe}.

\begin{thm}[\cite{GuMaPe}]\label{thm:gumape}
Let $\Sigma_{\varepsilon,\eta}(\Omega)$ be defined
by~\eqref{eq:cislo2}. Then there exist positive numbers
$\Delta_1,\Delta_2\in\Z[\eta]$ and a strictly increasing sequence
$(s_n)_{n\in\Z}$ such that
\begin{enumerate}
\item $\Sigma_{\varepsilon,\eta}(\Omega) = \{s_n \mid
n\in\Z\}\subset\Z[\eta]$.

\item $\Delta_1^*>0$, $\Delta_2^*<0$, $\Delta_1^*-\Delta_2^*\geq
l>\max(\Delta_1^*,-\Delta_2^*)$.

\item $s_{n+1}-s_n\in\{\Delta_1,\Delta_2,\Delta_1+\Delta_2\}$, for
all $n\in\Z$, and, moreover,
$$
s_{n+1}=\left\{\begin{array}{ll}
s_n+\Delta_1    &\hbox{ if } \quad s_n^*\in[c,c+l-\Delta_1^*)\,,\\[2mm]
s_n+\Delta_1+\Delta_2   &\hbox{ if } \quad s_n^*\in[c+l-\Delta_1^*,c-\Delta_2^*)\,,\\[2mm]
s_n+\Delta_2    &\hbox{ if } \quad s_n^*\in[c-\Delta_2^*,c+l)\,.
\end{array}\right.
$$

\item Numbers $\Delta_1$ and $\Delta_2$ depend only on parameters $\varepsilon,\eta$
and the length $l$ of the interval $\Omega$. In particular, they
do not depend on the position $c$ of $\Omega$ on the real line.

\end{enumerate}
\end{thm}

We see that the set $\{s_n^* \mid n\in\Z\}$ is an orbit under the
3iet with permutation $\pi=(321)$ and parameters $l-\Delta_1^*$,
$\Delta_1^*-\Delta_2^*-l$ and $l+\Delta_2^*$ (if
$l<\Delta_1^*-\Delta_2^*$) and it is an orbit under the 2iet with
permutation $\pi=(21)$ and parameters $l-\Delta_1^*$ and
$l+\Delta_2^*$ (if $l=\Delta_1^*-\Delta_2^*$). Thus every
cut-and-project sequence can be viewed as a geometric
representation of an orbit of a point under exchange of two or
three intervals.

The determination of $\Delta_1$, $\Delta_2$ is in general
laborious; the values $\Delta_1$, $\Delta_2$ depend on the
continued fraction expansions of parameters $\varepsilon$ of
$\eta$, according to the length $l$ of the acceptance window
$\Omega=[c,c+l)$.

In case that
\begin{equation}\label{eq:cislo3}
\varepsilon\in(0,1),\quad \eta>0 \quad\hbox{and}\quad 1\geq
l>\max(1-\varepsilon,\varepsilon)\,,
\end{equation}
one has
\begin{equation}\label{eq:cislo20}
\Delta_1=1+\eta \quad\hbox{and}\quad \Delta_2=\eta\,,
\end{equation}
i.e., the corresponding triple of shifts in the prescription of
the exchange of intervals is $\Delta_1^*=1-\varepsilon$,
$\Delta_1^*+\Delta_2^* = 1-2\varepsilon$,
$\Delta_2^*=-\varepsilon$. In fact, without loss of generality, we
can limit our consideration to cut-and-project sequences with
parameters satisfying~\eqref{eq:cislo3}, since in~\cite{GuMaPe} it
is shown that every cut-and-project sequence is equal to
$\mu\Sigma_{\varepsilon,\eta}(\Omega)$, where $\varepsilon$,
$\eta$ and length $l$ of the interval $\Omega$
satisfy~\eqref{eq:cislo3} and $\mu\in\R$. By that, we have shown
how to interpret a cut-and-project set as an orbit under an
exchange of 3 (or 2) intervals with the permutation (321) (or
(21)).

On the other hand, let us show that every exchange of three
intervals with permutation (321) can be represented geometrically
using a cut-and-project scheme. First realize that studying the
orbit of a point $x_0\in I$ under the 3iet $\widetilde{T}$
of~\eqref{eq:cislo4}, we can, without loss of generality,
substitute $\widetilde{T}$ by the transformation
$T(x)=\frac1{\mu}\widetilde{T}\bigl(\mu(x-c)\bigr)+c$ for
arbitrary $\mu,c\in\R$, $\mu\neq 0$, and instead of the orbit of
$x_0$ under $\widetilde{T}$ consider the orbit of the point
$y_0=c+\frac{x_0}{\mu}$ under the transformation $T$. In
particular, putting $\mu=\alpha_1+2\alpha_2+\alpha_3$ and
$c=-x_0\mu^{-1}$, we have the orbit of $y_0=0$ under the mapping
$T:[c,c+l)\mapsto[c,c+l)$
\begin{equation}\label{eq:cislo5}
T(x)=\left\{\begin{array}{ll}
 x+1-\varepsilon &\hbox{for }\ x\in[c,c+l-1+\varepsilon)\,,\\
 x+1-2\varepsilon &\hbox{for }\ x\in[c+l-1+\varepsilon,c+\varepsilon)\,,\\
 x-\varepsilon &\hbox{for }\ x\in[c+\varepsilon,c+l)\,,
\end{array}\right.
\end{equation}
where we have denoted by $\varepsilon$ and $l$ the new parameters
\begin{equation}\label{eq:cisloCCC}
\varepsilon:=\frac{\alpha_1+\alpha_2}{\alpha_1+2\alpha_2+\alpha_3}\quad\hbox{
and }\quad
l:=\frac{\alpha_1+\alpha_2+\alpha_3}{\alpha_1+2\alpha_2+\alpha_3}\,.
\end{equation}
Let us mention that under such parameters, the minimality property
of the transformation $T$ in~\eqref{eq:cislo5} is equivalent to
the requirement that $\varepsilon$ be irrational.

For the above defined values of $\varepsilon,l,c$ and arbitrary
irrational $\eta>0$ put $\Omega=[c,c+l)$ and consider the
cut-and-project set $\Sigma_{\varepsilon,\eta}(\Omega)$. Since
$0\in\Omega$, we have also
$0\in\Sigma_{\varepsilon,\eta}(\Omega)$. The strictly increasing
sequence $(s_n)_{n\in\Z}$ from Theorem~\ref{thm:gumape} can be
indexed in such a way that $s_0=0$. Since our parameters
$\varepsilon,l,\eta$ satisfy~\eqref{eq:cislo3} (and $l<1$), the
right neighbor $s_{n+1}$ of the point $s_n$ is given by the
position of $s_n^*$ in the interval $[c,c+l)$, namely by the
transformation $T(x)$. In particular, we have
$s_{n+1}^*=T(s_n^*)$. Therefore the set
$$
\{s_n^*\mid n\in\N\} =
\bigl(\Sigma_{\varepsilon,\eta}[c,c+l)\bigr)^* =
\Z[\varepsilon]\cap \Omega
$$
is the orbit of the point 0 under the transformation $T$.

Note that we have decided to consider instead of an orbit of an
arbitrary point under a 3iet $\widetilde{T}$ with the domain being
an interval starting at 0, the orbit of 0 under the 3iet $T$ given
by~\eqref{eq:cislo5}, with parameters $\varepsilon,l,c$ satisfying
\begin{equation}\label{eq:cislo6}
\varepsilon\in(0,1), \quad 1>l>\max(1-\varepsilon,\varepsilon),
\quad 0\in[c,c+l)\,.
\end{equation}
Let us summarize the advantages of such new notation:
\begin{itemize}

\item
Points of the sequence
$\bigl(T^n(0)\bigr)_{n\in\Z}\subset\Z[\varepsilon]$ which has a
chaotic behavior in the interval $[c,c+l)$ can be, using the star
map $*:\Z[\eta]\mapsto\Z[\varepsilon]$, represented by a strictly
increasing sequence $(s_n)_{n\in\Z}$ such that $s_n^*=T^n(0)$ for
all $n\in\Z$.

\item
The orbit of 0 can be simply expressed as
$$
\{T^n(0)\mid n\in\Z\} = \Z[\varepsilon] \cap [c,c+l)\,.
$$
For the orbit of an arbitrary point $x_0\in[c,c+l)$ under $T$, one
can write
\begin{equation}\label{eq:cislo32}
\{T^n(x_0)\mid n\in\Z\} = x_0 + \bigl(\Z[\varepsilon] \cap
[c-x_0,c+l-x_0)\bigr) = (x_0+\Z[\varepsilon])\cap [c,c+l)\,.
\end{equation}

\end{itemize}

Further advantages of the presented point of view on 3iets by
cut-and-project sequences will be clear from the following
section.

\begin{pozn}\label{pozn:jineEty}
To conclude the section, let us stress that for the 3iet $T$ the
parameter $\eta$ was chosen arbitrarily, except the requirement of
irrationality and positiveness. Then adjacency of points $x,y$,
$x<y$, in the set $\Sigma_{\varepsilon,\eta}(\Omega)$ indicates
that their star map images $x^*$, $y^*$ are consecutive iterations
of $T$, i.e., $T(x^*)=y^*$. Choosing the parameter $\eta<0$, we
obtain again a cut-and-project set
$\Sigma_{\varepsilon,\eta}(\Omega)$ but with different $\Delta_1$,
$\Delta_2$. Therefore the corresponding 3iet is different from
$T$. From the definition of a cut-and-project set, it can be
easily shown that
$$
\Sigma_{\varepsilon,\eta}(\Omega)
=\Sigma_{1-\varepsilon,1-\eta}(\Omega).
$$
Therefore in case that $\eta<-1$, the corresponding
cut-and-project set represents a 3iet, in which we interchange the
lengths of the first and last intervals, i.e., the mapping
$T^{-1}$. In fact, the `dangerous' choice for the irrational
parameter $\eta$ is $\eta\in(-1,0)$.
\end{pozn}

\section{First return map}

Let $T:I\mapsto I$ be a $k$-interval exchange transformation with
minimality property and let $J$ be an interval $J\subset I$, $J$
closed from the left and open from the right, say
$[\hat{c},\hat{c}+\hat{l})$.

The minimality property of $T$ ensures that for every $z\in J$
there exists a positive integer $i\in\N$ such that $T^{i}(z)\in
J$. The minimal such $i$ is called the return time of $z$ and
denoted by $r(z)$.

To every $z\in J$ we associate a `return name', i.e., a finite
word $w=v_0v_1\cdots v_{r(z)-1}$ in the alphabet $\{1,\dots,k\}$,
whose length is equal to the return time of $z$ and for all $i$,
$0\leq i<r(z)$ we have
$$
v_i=X \quad\hbox{if}\quad T^i(z)\in I_X\,.
$$
To the given subinterval $J$ of $I$, we define the map
$T_J:J\mapsto J$ by the prescription
$$
T_J(z) = T^{r(z)}(z)\,,
$$
which is called the first return map.

Since for a fixed interval $J$ the return time $r(z)$ is bounded,
there exist only finitely many return names. It is obvious, that
points $z\in J$ with the same return name form an interval, and
$J$ is thus a finite disjoint union of such subintervals, say
$J_1,\dots,J_p$. The boundary points of these intervals can be
easily described by the notion of ancestor in $J$.

The minimality property of $T$ ensures that for every $y\in I$
there exists $z\in J$ such that
$y\in\{z,T(z),\dots,T^{r(z)-1}(z)\}$. Such $z$ is uniquely
determined and we call it the ancestor of $y$ in the interval $J$.
We denote $z={\rm anc}_{\scriptscriptstyle J}(y)$.

The boundary points of the intervals $J_1,\dots, J_p$ are then
exactly the following points:
\begin{equation}\label{list}
\begin{array}{l}
 \bullet \quad \hat{c},\ \hat{c}+\hat{l}\ \hbox{(i.e., the boundary
points of $J$ itself);}\\
 \bullet \quad {\rm anc}_{\scriptscriptstyle J}(\hat{c}+\hat{l});\\
 \bullet \quad {\rm anc}_{\scriptscriptstyle J}(\alpha_1+\alpha_2+\cdots+\alpha_i)\ \hbox{
 for } i=1,2,\dots,k-1;\\
 \bullet \quad \hbox{and the point $z\in J$ such that }
 T^{r(z)}(z)=\hat{c}.
\end{array}
\end{equation}
This implies that for a $k$-iet the number of different return
names is at most $k+2$. It is obvious, that the first return map
$T_J$ is again a $m$-iet for some $m\leq k+2$. In fact, it is
known that $m\leq k+1$ (see~\cite{CFS}, Chap.~5). For a 3iet which
we study in this paper, we can say even more. The following
theorem is a direct consequence of Theorem~\ref{thm:gumape}.

\begin{thm}
Let $T:I\mapsto I$ be a 3iet with permutation (321) and satisfying
minimality property, and let $J\subset I$ be an interval. Then the
first return map $T_J$ is either a 3iet with permutation (321) or
a 2iet with permutation (21).
\end{thm}

\section{First return map and substitution invariance}

Let us now see how the notions of first return map, return time
and return name are related to substitution invariance of words
coding 3iet. We will focus on non-degenerate 3iet words. Let us
mention that non-degeneracy in terms of parameters $\varepsilon,l$
of~\eqref{eq:cisloCCC} means that $l\notin\Z[\varepsilon]$,
cf.~\eqref{eq:cisloCC}.

Consider a 3iet $T:[c,c+l)\mapsto[c,c+l)$ of~\eqref{eq:cislo5}
with parameters~\eqref{eq:cislo6} and an interval
$J\subset[c,c+l)$ such that $0\in J$. Let $w_1,\dots,w_p$ be all
possible return names of points $z\in J$. Then the infinite word
$u=(u_n)_{n\in\Z}$ coding 0 under the transformation $T$ can be
written as a concatenation
\begin{equation}\label{eq:cislo12}
u=\cdots w_{j_{-2}}w_{j_{-1}}|w_{j_0}w_{j_1}w_{j_2}\cdots\,,\qquad
\hbox{ with } \ j_i\in\{1,\dots,p\}\,.
\end{equation}
The starting letters of the blocks $w_{j_m}$ correspond to
positions $n$ in the infinite word $u$ if and only if $T^n(0)\in
J$. More formally, we have
$$
w_{j_m}w_{j_{m+1}}w_{j_{m+2}}\cdots =
u_nu_{n+1}u_{n+2}u_{n+3}\cdots \qquad\iff\qquad T^n(0)\in J\,.
$$
Suppose we have an interval $J\subset I$, $0\in J$ such that the
first return map $T_J$ satisfies
\begin{enumerate}
\item[P1.] $T_J$ is homothetic with $T$, i.e.,
$$
T_J(x)=\nu T(\tfrac{x}{\nu})\,,\quad\hbox{for }\ x\in J \hbox{ and
some }\ \nu\in(-1,1)\,,
$$
which means that $T_J$ is an exchange of intervals $J_1=\nu I_1$,
$J_2=\nu I_2$, and $J_3=\nu I_3$;
\item[P2.] the set of return names defined by $J$ has three elements.
\end{enumerate}
Then the sequence of indices $(j_m)_{m\in\Z}$ defining the
ordering of finite words $w_1,w_2,w_3$ in the
concatenation~\eqref{eq:cislo12} equals to the infinite word $u$.
In particular, it means that $u$ is invariant under the
substitution
$$
\begin{array}{ccl}
1&\mapsto&\varphi(1)=w_1,\\
2&\mapsto&\varphi(2)=w_2,\\
3&\mapsto&\varphi(3)=w_3.
\end{array}
$$

We stand therefore in front of the following questions: How to
decide, for which 3iets a subinterval $J\subset I$ with properties
P1.\ and P2.\ exists? What can be said in case that such $J$ does
not exist?

In case that $u=(u_n)_{n\in\Z}$ is a non-degenerate 3iet word
coding the orbit of 0 under the transformation $T$ defined
by~\eqref{eq:cislo5}, the second question is solved by the
paper~\cite{bluetheorem}, as follows.

The existence of a substitution $\varphi$ over the alphabet
$\{1,2,3\}$, under which the word $u$ is invariant, means that $u$
can be written as a concatenation of blocks $\varphi(1)$,
$\varphi(2)$, $\varphi(3)$, i.e.,
\begin{equation}\label{eq:cislo8}
u\quad=\quad\cdots u_{-2}u_{-1}|u_0u_1u_2\cdots \quad=\quad \cdots
\varphi(u_{-2})\varphi(u_{-1})\mid\varphi(u_0)\varphi(u_1)\varphi(u_2)\cdots\,.
\end{equation}
In~\cite{bluetheorem} one considers a non-degenerate 3iet word $u$
invariant under a primitive substitution $\varphi$ and studies for
$i=1,2,3$ the set $E_{\varphi(i)}$ of points $T^n(0)$ such that
the block $\varphi(i)$ starts at position $n$ in the
concatenation~\eqref{eq:cislo8}. Formally,
$$
E_{\varphi(i)}= \{T^n(0) \mid \exists m\in\Z,\
\varphi(i)\varphi(u_m)\varphi(u_{m+1})\cdots =
u_nu_{n+1}u_{n+2}\cdots\}\,.
$$
As a result, several properties of a matrix of substitution
$\varphi$ are described. Recall that for a substitution $\varphi$
over the alphabet $\A=\{1,2,\dots, k\}$ one defines the
substitution matrix $M_\varphi$ by
$$
(M_\varphi)_{ij} = \hbox{ number of letters $i$ in the word
$\varphi(j)$}\,,\quad 1\leq i,j \leq k\,.
$$
Such matrix has obviously non-negative integer entries and if the
substitution $\varphi$ is primitive, the matrix $M_\varphi$ is
primitive as well, and therefore one can apply the
Perron-Frobenius theorem.

 We summarize several statements
of~\cite{bluetheorem} in the following theorem.

\begin{thm}[\cite{bluetheorem}]\label{thm:ABMP}
Let $u=(u_n)_{n\in\Z}$ be a non-degenerate 3iet word with
parameters $\varepsilon,l,c$ satisfying~\eqref{eq:cislo6}. Let
$\varphi$ be a primitive substitution such that $\varphi(u)=u$.
Then
\begin{enumerate}
\item[(i)] $\varepsilon$ is a Sturm number, i.e., $\varepsilon$ is
a quadratic irrational in $(0,1)$ such that its algebraic
conjugate $\varepsilon'$ satisfies $\varepsilon'\notin(0,1)$;

\item[(ii)] the dominant eigenvalue $\Lambda$ of the
matrix $M_\varphi$ of the substitution $\varphi$ is a quadratic
unit in $\Q(\varepsilon)$;

\item[(iii)] the column vector
$(1-\varepsilon,1-2\varepsilon,-\varepsilon)^T$ is a right
eigenvector of $M_\varphi$ corresponding to the eigenvalue
$\Lambda'$, i.e., to the algebraic conjugate of $\Lambda$;

\item[(iv)] parameters $c,l\in\Q(\varepsilon)$;

\item[(v)] $E_{\varphi(i)} = \Lambda'\bigl(I_i\cap \Z[\varepsilon]\bigr)$ for
$i=1,2,3$.
\end{enumerate}
\end{thm}

The statement (v) in particular says that the existence of a
substitution $\varphi$ under which a non-degenerate 3iet word $u$
is invariant forces existence of an interval $J\subset I$ with
properties P1.\ and P2. We have already explained that existence
of an interval $J$ with properties P1.\ and P2.\ forces
substitution invariance. We have thus the following statement.

\begin{prop}\label{p:ekvival}
Let $u=(u_n)_{n\in\Z}$ be a non-degenerate 3iet word with
parameters $\varepsilon,l,c$ satisfying~\eqref{eq:cislo6}. Then
there exists a primitive substitution $\varphi$ under which $u$ is
invariant, if and only if there exists an interval $J\subset I$
with properties P1.\ and P2.
\end{prop}

Let us first derive two simple observations which complement
results of~\cite{bluetheorem}.

\begin{lem}\label{l:1}
For $\Lambda,\Lambda'$ and $\varepsilon$ from
Theorem~\ref{thm:ABMP} we have
$$
\Lambda\Z[\varepsilon] = \Lambda'\Z[\varepsilon] =
\Z[\varepsilon]\,.
$$
\end{lem}

\pf Statement (iii) of Theorem~\ref{thm:ABMP} implies
$$
M_\varphi\left(\!\!\!\begin{array}{c}1-\varepsilon\\
1-2\varepsilon\\-\varepsilon\end{array}\!\!\!\right)
 = \Lambda'\left(\!\!\!\begin{array}{c}1-\varepsilon\\
1-2\varepsilon\\-\varepsilon\end{array}\!\!\!\right)\,.
$$
Since $M_\varphi$ is an integer matrix, we obtain from the third
row of the above equality that
$\Lambda'\varepsilon\in\Z[\varepsilon]$. Subtracting third row
from the first one we get $\Lambda'\in\Z[\varepsilon]$. Since
$\Z[\varepsilon]$ is closed under addition, we have
$\Lambda'\Z[\varepsilon]\subseteq\Z[\varepsilon]$.

Since $\Lambda$ is a quadratic integer, we have
$\Lambda+\Lambda'\in\Z$. This implies that
$\Lambda\in\Z-\Lambda'\in\Z[\varepsilon]$, whence
$\Lambda\varepsilon\in\varepsilon\Z-\Lambda'\varepsilon\in\Z[\varepsilon]$,
and thus $\Lambda\Z[\varepsilon]\subseteq\Z[\varepsilon]$.

Now since $\Lambda$ is a unit, we have $\Lambda\Lambda'=\pm1$, and
therefore multiplying
$\Lambda\Z[\varepsilon]\subseteq\Z[\varepsilon]$ by $\Lambda'$ we
obtain $\Z[\varepsilon]\subseteq\Lambda'\Z[\varepsilon]$. \pfk

It is obvious that in our considerations, $\varepsilon$ must be a
quadratic irrational. When putting a 3iet with such a parameter
into context of cut-and-project sets, we need to specify the slope
of the second projection, i.e., the parameter $\eta$. Choosing
$\eta=-\varepsilon'$, where $\varepsilon'$ is the algebraic
conjugate of $\varepsilon$,  the star map $x=a+b\eta\mapsto
x^*=a-b\varepsilon$ becomes the Galois automorphism in
$\Q(\varepsilon)$. We will use the notation $x=a+b\varepsilon$,
$a,b\in\Q$ $\mapsto x'=a+b\varepsilon'$, as is usual. Recall that
for $x,y\in\Q(\varepsilon)$ we have
$$
(x+y)'=x'+y' \quad\hbox{ and }\quad (xy)'=x'y'\,.
$$
With such notation, $\Sigma_{\varepsilon,-\varepsilon'}(\Omega)$
can be rewritten in the form
\begin{equation}\label{eq:cislo70}
\Sigma_{\varepsilon,-\varepsilon'}(\Omega) =
\{x\in\Z[\varepsilon']\mid x'\in\Omega\} \,.
\end{equation}

\begin{lem}\label{l:2}
Let $\varepsilon$ be a quadratic irrational and let $\Lambda$ be a
quadratic unit in $\Q(\varepsilon)$ such that
\begin{equation}\label{eq:cislo10}
\Lambda\Z[\varepsilon]=\Z[\varepsilon]\,.
\end{equation}
\begin{itemize}
\item
Then for any acceptance window $\Omega$ we have
$$
\Lambda\Sigma_{\varepsilon,-\varepsilon'}(\Omega) =
\Sigma_{\varepsilon,-\varepsilon'}(\Lambda'\Omega)\,.
$$
\item
If moreover $\varepsilon'<0$, $\Lambda>1$, $\Lambda'\in(0,1)$ and
$T:[c,c+l)\mapsto[c,c+l)$ is a 3iet with parameters
satisfying~\eqref{eq:cislo6}, then the first return map $T_J$ for
the interval $J=\Lambda'[c,c+l)$ is a 3iet homothetic with $T$.
\end{itemize}
\end{lem}

\pf Since $\Lambda\Lambda'=\pm1$, multiplying
of~\eqref{eq:cislo10} by $\Lambda'$ leads to
$\Lambda'\Z[\varepsilon]=\Z[\varepsilon]=\Z[-\varepsilon]$. By
algebraic conjugation we obtain
$\Lambda\Z[\varepsilon']=\Z[\varepsilon']=\Z[-\varepsilon']$. Note
that in general $\Z[\varepsilon]\neq\Z[\varepsilon']$.
From~\eqref{eq:cislo70} we obtain
$$
\begin{aligned}
\Lambda\Sigma_{\varepsilon,-\varepsilon'}(\Omega) =
\Lambda\{x\in\Z[\varepsilon']\mid x'\in\Omega\} &= \{\Lambda
x\in\Z[\varepsilon']\mid \Lambda'x'\in\Lambda'\Omega\} =\\&=
\{y\in\Z[\varepsilon']\mid y'\in\Lambda'\Omega\} =
\Sigma_{\varepsilon,-\varepsilon'}(\Lambda'\Omega).
\end{aligned}
$$

This however means that the distances between adjacent elements of
the cut-and-project set
$\Sigma_{\varepsilon,-\varepsilon'}(\Lambda'\Omega)$ are $\Lambda$
multiples of the distances between adjacent elements of the
cut-and-project set $\Sigma_{\varepsilon,-\varepsilon'}(\Omega)$.
Since the star map images (in our case the images under the Galois
automorphism) of the distances between neighbors in a
cut-and-project set correspond to translations in the
corresponding 3iet (see Theorem~\ref{thm:gumape}), the factor of
homothety between the two 3iets is $\Lambda$.

If the parameter $\eta=-\varepsilon'>0$, the  3iet mappings
corresponding to $\Sigma_{\varepsilon,-\varepsilon'}(\Omega)$ and
$\Sigma_{\varepsilon,-\varepsilon'}(\Lambda'\Omega)$ are precisely
$T$ and $T_J$ respectively, see Remark~\ref{pozn:jineEty}.
 \pfk

Using Lemma~\ref{l:1} and statement (v) of
Theorem~\ref{thm:ABMP},we obtain
\begin{equation}\label{eq:cislo7}
E_{\varphi(i)} = (\Lambda'I_i) \cap \Z[\varepsilon] = (\Lambda'
I_i) \cap \{T^n(0) \mid n\in\Z\}\,.
\end{equation}

We are now in position to prove the main theorem of this section,
which provides a necessary and sufficient condition for
substitution invariance of a non-degenerate 3iet word.

\begin{prop}\label{thm:P}
Let $u$ be a non-degenerate 3iet word with parameters
$\varepsilon,l,c$, such that $\varepsilon$ is a Sturm number
having $\varepsilon'<0$ and $l,c\in\Q(\varepsilon)$,
$l\notin\Z[\varepsilon]$. Then $u$ is invariant under a primitive
substitution if and only if there exists a quadratic unit
$\Lambda\in\Q(\varepsilon)$, $\Lambda>1$, with conjugate
$\Lambda'\in(0,1)$, such that
 \begin{enumerate}
 \item[C1.] $\Lambda\Z[\varepsilon]=\Z[\varepsilon]$, and
 \item[C2.] for the interval $J=\Lambda'[c,c+l)$, one has
 $$
 {\rm anc}_{\scriptscriptstyle J}(c+\varepsilon),{\rm anc}_{\scriptscriptstyle J}(c+l-(1-\varepsilon))\in
 \bigl\{\Lambda'c,\Lambda'(c+\varepsilon),\Lambda'(c+l-(1-\varepsilon))\bigr\}.
 $$
 \end{enumerate}
\end{prop}

\pf Let $u$ be invariant under a primitive substitution $\varphi$.
We search for $\Lambda$ with properties C1.\ and C2.\ of the
proposition. According to Theorem~\ref{thm:ABMP}, the dominant
eigenvalue of the matrix $M_\varphi$ is a quadratic unit in
$\Q(\varepsilon)$, i.e., its conjugate belongs to the interval
$(-1,1)$. If the conjugate is positive, we use for $\Lambda$ the
dominant eigenvalue of $M_\varphi$. Otherwise, since $u$ is
invariant also under the substitution $\varphi^2$, we take for
$\Lambda$ the dominant eigenvalue of the matrix
$M_{\varphi^2}=M_{\varphi}^2$.

The validity of property C1.\ follows from Lemma~\ref{l:1}.
Equation~\eqref{eq:cislo7} states that the interval $J=\Lambda'I$
defines only three return names and that the subintervals
corresponding to these return names are $\Lambda'I_1$,
$\Lambda'I_2$ and $\Lambda'I_3$. Since $I=[c,c+l)$, these are
$\Lambda'I_1=\bigl[\Lambda'c,\Lambda'(c+l-1+\varepsilon)\bigr)$,
$\Lambda'I_2=\bigl[\Lambda'(c+l-1+\varepsilon),\Lambda'(c+\varepsilon)\bigr)$,
and
$\Lambda'I_1=\bigl[\Lambda'(c+\varepsilon),\Lambda'(c+l)\bigr)$.
The list~\eqref{list} defines the boundary points of subintervals
determining the return names. Property C2.\ follows.

For the opposite implication, realize that by Lemma~\ref{l:2}
property C1.\ ensures that $T_J$ is a 3iet with subintervals
$\Lambda'[c,c+l-1+\varepsilon)$
$\Lambda'[c+l-1+\varepsilon,c+\varepsilon)$, and
$\Lambda'[c+\varepsilon, c+l)$. This, together with property C2.,
forces that points of the list~\eqref{list} belong to the set
$\{\Lambda'c,\Lambda'(c+\varepsilon),\Lambda'(c+l-1+\varepsilon)\}$,
and thus the interval $J=\Lambda'I$ defines three return names.
Hence according to Proposition~\ref{p:ekvival}, the infinite word
$u$ is invariant under a primitive substitution. \pfk

\begin{pozn}\label{pozn:J}
The proof of the above proposition directly implies that in case
that $u$ is invariant under a substitution $\varphi$, the scaling
factor $\Lambda$ from Proposition~\ref{thm:P} can be taken to be
the dominant eigenvalue of the substitution matrix $M_\varphi$ or
$M_{\varphi^2}=M_\varphi^2$.
\end{pozn}

\section{Characterization of substitution invariant 3iet words}

We now have to solve the question, when for a given Sturm number
$\varepsilon$ and parameters $c,l\in\Q(\varepsilon)$
satisfying~\eqref{eq:cislo6} there exists $\Lambda$ with
properties C1.\ and C2.\ of Proposition~\ref{thm:P}. Finding
$\Lambda$ having the first of the properties is simple.

\begin{lem}\label{l:3}
Let $\varepsilon$ be irrational, solution of the equation
$Ax^2+Bx+C=0$. Then there exists a quadratic unit
$\Lambda\in\Q(\varepsilon)$ such that
\begin{equation}\label{eq:cislo33}
\Lambda>1, \ \Lambda'\in(0,1),\ \hbox{and} \
\Lambda\Z[\varepsilon]=\Z[\varepsilon]\,.
\end{equation}
\end{lem}

\pf Let the pair of integers $X,Y$ be a non-trivial solution of
the Pell equation
$$
X^2-(B^2-4AC)Y^2=1\,.
$$
Put $\gamma:=X+BY+2AY\varepsilon$. Using
$A\varepsilon^2=-B\varepsilon-C$, we easily verify that
$\gamma\varepsilon\in\Z[\varepsilon]$. Using
$A(\varepsilon+\varepsilon')=-B$ and $A\varepsilon\varepsilon'=C$,
we derive that $\gamma\gamma'=1$. This implies
$$
\gamma\Z[\varepsilon]=\gamma'\Z[\varepsilon]=\Z[\varepsilon].
$$
Finally, we put $\Lambda=\max\{|\gamma|,|\gamma'|\}$.
 \pfk

In Lemma~\ref{l:3} we have found $\Lambda$ with property C1. It is
more difficult to decide when $\Lambda$ satisfies also property
C2.\ of Proposition~\ref{thm:P}. By definition of the map $T$, it
follows that $x$ and $T(x)$ differ by an element of
$\Z[\varepsilon]$. Therefore for arbitrary $z_0$ and its ancestor
${\rm anc}_{\scriptscriptstyle J}(z_0)$ we have $z_0-{\rm
anc}_{\scriptscriptstyle J}(z_0)\in\Z[\varepsilon]$. It is useful
to introduce an equivalence on $\Q(\varepsilon)$ as follows. We
say that elements $x,y\in\Q(\varepsilon)$ are equivalent if their
difference belongs to $\Z[\varepsilon]$. Formally,
$$
x-y\in\Z[\varepsilon] \qquad\iff\qquad x\sim y\,.
$$

For the parameters $c,l\in\Q(\varepsilon)$, one can find $q\in\N$
such that $c,l\in\frac{1}{q}\Z[\varepsilon]$. Clearly, ${\rm
anc}_{\scriptscriptstyle J}(c+\varepsilon)$ and ${\rm
anc}_{\scriptscriptstyle J}(c+l-1+\varepsilon)$ also belong to the
set $\frac{1}{q}\Z[\varepsilon]$. The set to which belong
ancestors of $c+\varepsilon$ and $c+l-1+\varepsilon$ can be
restricted even more. For, the equivalence $\sim$ divides the set
$\frac{1}{q}\Z[\varepsilon]$ into $q^2$ classes of equivalence of
the form
$$
T_{ij}:=\frac{i+j\varepsilon}{q}+\Z[\varepsilon]\,,\quad \hbox{
where } 0\leq i,j\leq q-1\,.
$$
Relation $\Lambda'\Z[\varepsilon]=\Z[\varepsilon]$ implies
$$
z\in\Z[\varepsilon] \quad\iff\quad \Lambda'z\in\Z[\varepsilon]\,.
$$
Therefore the mapping $\psi(T_{ij})=\Lambda'T_{ij}$ is a bijection
on the set of $q^2$ classes of equivalence. For every bijection
$\psi$ on a finite set, there exists an iteration $s\in\N$, $s\geq
1$, such that $\psi^s={\rm id}$. Denoting $L:=\Lambda^s$, the
number $L$ has obviously all properties of $\Lambda$, namely
\begin{itemize}
\item[a)] $L$ is a quadratic unit in $\Q(\varepsilon)$;
\item[b)] $L>1$, $L'\in(0,1)$;
\item[c)] $L\Z[\varepsilon]=\Z[\varepsilon]$;
\end{itemize}
and moreover
\begin{itemize}
\item[d)] $L'\bigl(\frac{i+j\varepsilon}{q}+\Z[\varepsilon]\bigr)=
\frac{i+j\varepsilon}{q}+\Z[\varepsilon]$, \quad for all $i,j$,
$1\leq i,j\leq q-1$.
\end{itemize}

Having a quadratic unit $\Lambda$ with properties of the number
$L$ in items a) -- d), it is less difficult to decide about
validity of the condition
\begin{equation}\label{eq:cislo38}
{\rm anc}_{\scriptscriptstyle J}(c+\varepsilon), {\rm
anc}_{\scriptscriptstyle J}(c+l-1+\varepsilon) \ \in \
\bigl\{\Lambda'c,\Lambda'(c+\varepsilon),
\Lambda'(c+l-1+\varepsilon)\bigr\}\,.
\end{equation}
Non-degeneracy of the infinite word $u$ implies that
$l\notin\Z[\varepsilon]$, and therefore $c+\varepsilon\ \not\sim\
c+l-1+\varepsilon$. Since for every $z_0\in\frac1q\Z[\varepsilon]$
we have now
$$
z_0\ \sim\ {\rm anc}_{\scriptscriptstyle J}(z_0) \ \sim \
\Lambda'z_0\,,
$$
the condition~\eqref{eq:cislo38} in fact means
\begin{equation}\label{eq:cislo39}
{\rm anc}_{\scriptscriptstyle J}(c+l-1+\varepsilon) =
\Lambda'(c+l-1+\varepsilon)
\end{equation}
and
\begin{equation}\label{eq:cislo40}
{\rm anc}_{\scriptscriptstyle J}(c+\varepsilon) \ \in \
\bigl\{\Lambda'c,\Lambda'(c+\varepsilon)\bigr\}\,.
\end{equation}

\begin{lem}\label{l:H}
Let $\varepsilon$ be a Sturm number with $\varepsilon'<0$. Let
$l,c\in\frac1q\Z[\varepsilon]$. Let $\Lambda$ satisfy properties
of $L$ in a) -- d) and let $J=\Lambda'[c,c+l)$. Then for arbitrary
$z_0\in\frac1q\Z[\varepsilon]\cap[c,c+l)$, one has
$$
{\rm anc}_{\scriptscriptstyle J}(z_0)=\Lambda'z_0 \quad\iff\quad
z'_0\leq 0\leq (T(z_0))'\,.
$$
\end{lem}

\pf The transformation $T$ preserves the classes of equivalence
and thus for the orbit of a point $z_0$ it holds that
$$
\{T^n(z_0)\mid n\in\Z\} \ \subset \ z_0+\Z[\varepsilon]\,.
$$
As
$(T^{n+1}(z_0)-T^n(z_0))'\in\{1-\varepsilon',1-2\varepsilon',-\varepsilon'\}$,
the assumption $\varepsilon'<0$ implies that the sequence
$(s_n)_{n\in\Z}$,
$$
s_n:=(T^n(z_0))'
$$
is strictly increasing. By~\eqref{eq:cislo32} we have moreover
$$
\{T^n(z_0)\mid n\in\Z\} \ = \ \{s'_n\mid n\in\Z\} \ = \
(z_0+\Z[\varepsilon])\cap[c,c+l)\,.
$$
Since $0\in[c,c+l)$ and $\Lambda'\in(0,1)$, it is
$\Lambda'[c,c+l)\subset[c,c+l)$. This inclusion together with
property d) implies
$$
\{s'_n\mid n\in\Z\} \ \supset \ \Lambda'
\Bigl((z_0+\Z[\varepsilon])\cap[c,c+l)\Bigr) \ = \
\{\Lambda's'_n\mid n\in\Z\}\,.
$$
The strictly increasing sequence $(\Lambda s_n)_{n\in\Z}$ is
therefore a subsequence of the strictly increasing sequence
$(s_n)_{n\in\Z}$. Thus there exists a unique index $m$ such that
\begin{equation}\label{eq:cislo40bis}
\Lambda s_m \leq s_0 < s_1 \leq \Lambda s_{m+1}\,.
\end{equation}
For determination of the ancestor of the point $z_0=s'_0$ by
definition, we search for the maximal non-positive index $k\in\Z$
such that $T^k(z_0)\in\Lambda'[c,c+l)$, i.e., such that $T^k(z_0)$
is an element of the sequence $(\Lambda's'_n)_{n\in\Z}$. Since
both $(s_n)_{n\in\Z}$ and $(\Lambda s_n)_{n\in\Z}$ are strictly
increasing, we have $(T^k(z_0))' = \Lambda s_m$ and thus ${\rm
anc}_{\scriptscriptstyle J}(s'_0)=\Lambda's'_m$. Denoting
$s'_m=y_0$, equation~\eqref{eq:cislo40bis} can be rewritten
\begin{equation}\label{eq:cislo41}
\Lambda y'_0 \leq z'_0 < (T(z_0))' \leq \Lambda (T(y_0))'\,.
\end{equation}
On the other hand, recall that $\{s'_n\mid n\in\Z\} =
(z_0+\Z[\varepsilon])\cap[c,c+l)$ and the index $m$ for
which~\eqref{eq:cislo40bis} holds, is determined uniquely.
Therefore we can claim that ${\rm anc}_{\scriptscriptstyle
J}(z_0)=\Lambda'y_0$ if and only if $y_0$ verifies
inequalities~\eqref{eq:cislo41}. Thus ${\rm
anc}_{\scriptscriptstyle J}(z_0)=\Lambda'z_0$ if and only if
\begin{equation}\label{eq:cislo412}
\Lambda z'_0 \leq z'_0 < (T(z_0))' \leq \Lambda (T(z_0))'\,.
\end{equation}
Note that strict inequality in the middle is trivial and it is
satisfied by arbitrary $z_0$. Since $\Lambda>1$ we have $\Lambda
z'_0\leq z'_0 \Leftrightarrow z'_0\leq 0$ and $(T(z_0))' \leq
\Lambda (T(z_0))' \Leftrightarrow 0\leq (T(z_0))'$, which
completes the proof.
 \pfk

\begin{thm}\label{t:hlavni}
Let $u$ be a non-degenerate 3iet word coding the orbit of the
point $x_0$ under a 3iet with permutation (321) and parameters
$\alpha_1,\alpha_2,\alpha_3$. Put
$$
\varepsilon:=\frac{\alpha_1+\alpha_2}{\alpha_1+2\alpha_2+\alpha_3}\,,\quad
l:=\frac{\alpha_1+\alpha_2+\alpha_3}{\alpha_1+2\alpha_2+\alpha_3}\,,
\quad\hbox{ and }\quad
c:=\frac{-x_0}{\alpha_1+2\alpha_2+\alpha_3}\,.
$$
Then $u$ is invariant under a primitive substitution if and only
if
\begin{enumerate}
\item $\varepsilon$ is a Sturm number;
\item $c,l\in\Q(\varepsilon)$;
\item $
\min(\varepsilon',1-\varepsilon')\ \leq \ c', \ c'+l' \ \leq \
\max(\varepsilon',1-\varepsilon')$.
%
\end{enumerate}
\end{thm}

\pf Theorem~\ref{thm:ABMP} claims that items 1.\ and 2.\ are
necessary conditions for existence of a primitive substitution
under which $u$ be invariant. Therefore we shall prove the
following statement:

If $\varepsilon$ is a Sturm number and $c,l\in\Q(\varepsilon)$,
then $u$ is invariant under a primitive substitution if and only
if condition 3.\ holds.

Note that the infinite word
$$
\cdots u_{-3}u_{-2}u_{-1}|u_0u_1u_2\cdots
$$
is substitution invariant if and only if
$$
\cdots u_{2}u_{1}u_{0}|u_{-1}u_{-2}u_{-3}\cdots
$$
is substitution invariant. At the same time, $\cdots
u_{2}u_{1}u_{0}|u_{-1}u_{-2}u_{-3}\cdots$ is a 3iet word coding
the transformation $T^{-1}$, i.e.\ the 3iet with parameters
$1-\varepsilon,l,c$. The fact that $\varepsilon$ is a Sturm number
means either $\varepsilon'<0$ or $\varepsilon'>1$. Instead of
parameters $\varepsilon,l,c$ we can thus have $1-\varepsilon,l,c$,
and therefore limit our study (without loss of generality) to
Sturm number $\varepsilon$ satisfying $\varepsilon'<0$. In that
case, inequalities in item 3.\ of the theorem are of the form
\begin{equation}\label{eq:cislo50}
\varepsilon' \leq c'+l' \leq 1-\varepsilon' \quad\hbox{and}\quad
\varepsilon' \leq -c' \leq 1-\varepsilon'\,.
\end{equation}

For the implication $\Rightarrow$, suppose that $u$ is invariant
under a primitive substitution $\varphi$. Denote $q\in\N$, such
that $c,l\in\frac1q\Z[\varepsilon]$. Since $u$ is invariant under
an arbitrary power of the substitution $\varphi$, we can use
Proposition~\ref{thm:P} and Remark~\ref{pozn:J} to find a number
$\Lambda$ with properties of the number $L$ described in a) -- d),
and such that for the interval $J=\Lambda'[c,c+l)$
equalities~\eqref{eq:cislo39} and~\eqref{eq:cislo40} hold.

When applying Lemma~\ref{l:H} on $z_0=c+l-1+\varepsilon$,
equality~\eqref{eq:cislo39} is equivalent to
$$
c'+l'-1+\varepsilon' \leq 0 \leq \bigl(T(c+l-1+\varepsilon)\bigr)'
= c'+l'-\varepsilon'\,,
$$
which is one of the inequalities in~\eqref{eq:cislo50}.

Now, let us study validity of~\eqref{eq:cislo40}. Since
$T(c+\varepsilon)=c$, we have ${\rm anc}_{\scriptscriptstyle
J}(c)={\rm anc}_{\scriptscriptstyle J}(c+\varepsilon)$.
Relation~\eqref{eq:cislo40} states that ${\rm
anc}_{\scriptscriptstyle J}(z_0)=\Lambda'z_0$ holds either for
$z_0=c$ or for $z_0=c+\varepsilon$. Therefore we have
\begin{equation}\label{eq:cislo43}
c'+\varepsilon' \leq 0\leq (T(c+\varepsilon))' \quad\hbox{or}\quad
c' \leq 0 \leq (T(c))'\,.
\end{equation}
Since $(T(c+\varepsilon))'=c'$ and $(T(c))'=c'+1-\varepsilon'$,
verifying at least one of the inequalities in~\eqref{eq:cislo43}
means
$$
c'+\varepsilon'\leq 0\leq c'+1-\varepsilon'\,,
$$
which is the other inequality in~\eqref{eq:cislo50}.

In the opposite implication, we take $\varepsilon$ Sturm and
$c,l\in\frac1q\Z[\varepsilon]$ and with the use of Lemma~\ref{l:3}
we find $\Lambda$ with properties of $L$ given in a) -- d). By
Lemma~\ref{l:H}, validity of inequalities~\eqref{eq:cislo50} is
equivalent to validity of~\eqref{eq:cislo39}
and~\eqref{eq:cislo40}. Therefore using Proposition~\ref{thm:P},
$u$ is invariant under a primitive substitution.
 \pfk

\begin{pozn}
Note that in the proof of the theorem we have applied
Lemma~\ref{l:H} only to points $z_0=c+l-1+\varepsilon$ and
$z_0=c+\varepsilon$. Realize that in fact, we do not need that
$\Lambda$ satisfies property d) of $L$, but only that
\begin{itemize}
\item[d')] $\Lambda'\bigl(c+\Z[\varepsilon]\bigr)=
c+\Z[\varepsilon]$ \quad and \quad
$\Lambda'\bigl(c+l+\Z[\varepsilon]\bigr)= c+l+\Z[\varepsilon]$.
\end{itemize}
This can be important when we search for minimal $\Lambda>1$ with
desired properties.
\end{pozn}

\section{Characterization of substitution invariant 3iet words using Sturmian words}

Comparing  Theorems   \ref{t:hlavni} and \ref{t:sturmian} we
immediately see a striking narrow connection between 3iet words
and Sturmian words, namely that the 3iet word $u=(u_n)_{n\in\N}$
is invariant under a primitive substitution if and only if  the
Sturmian  word with slope $\varepsilon$ and intercept $-c$ and the
Sturmian  word with slope $\varepsilon$ and intercept $\ell+c$ are
both invariant under a substitution.

In fact, as shown in~\cite{3iet}, these two Sturmian words appear
naturally as images of the given 3iet word by the following
morphisms.

\noindent Let us denote by
$\sigma_{01}:\{A,B,C\}^*\rightarrow\{0,1\}^*$ the
 morphism given by
 \begin{equation}\label{eq:delta}
 A \mapsto 0\,, \qquad
 B \mapsto 01\,, \qquad
 C \mapsto 1\,,
 \end{equation}
and by $\sigma_{10}:\{A,B,C\}^*\rightarrow\{0,1\}^*$ the  morphism
given by
 \begin{equation}
 A \mapsto 0\,, \qquad
 B \mapsto 10\,, \qquad
 C \mapsto 1\,.
 \end{equation}
One can verify that if  $u = (u_n)_{n\in\N}$ is a non-degenerate
3iet word with parameters $\varepsilon, \ell,c$
satisfying~\eqref{eq:cislo6}, then the infinite word
$$
\sigma_{01}(u) =
\sigma_{01}(u_0)\sigma_{01}(u_1)\sigma_{01}(u_2)\dots
$$
is the Sturmian word with slope $\varepsilon$ and intercept $-c$
and the infinite word $\sigma_{10}(u)$  is the Sturmian word with
slope $1-\varepsilon$ and intercept $\ell+c$.

With the definition of morphisms $\sigma_{01}$ and $\sigma_{10}$,
we can give the characterization of substitution invariant
non-degenerate 3iet words without use of any parameters.

\begin{coro} Let $u = (u_n)_{n\in\N}$ be a non-degenerate 3iet word.
Then $u$ is invariant under a primitive substitution  if and only
if both Sturmian words $\sigma_{10}(u)$ and $\sigma_{01}(u)$ are
invariant under substitution.
\end{coro}


\section*{Acknowledgements}

The authors acknowledge financial support by Czech Science
Foundation GA \v{C}R 201/05/0169, and by the grant LC06002 of the
Ministry of Education, Youth, and Sports of the Czech Republic.



\end{document}